\documentclass[11pt]{amsart}
\theoremstyle{plain}
\newtheorem{theorem}                 {THEOREM}
\newtheorem{proposition}  [theorem]  {PROPOSITION}

\newtheorem{lemma}        [theorem]  {LEMMA}

\theoremstyle{definition}

\title[COTANGENT BUNDLES WITH GENERAL NATURAL K\"AHLER STRUCTURES]{\bf{COTANGENT BUNDLES WITH GENERAL NATURAL K\"AHLER STRUCTURES}}

\author{S.~L.~Dru\c t\u a}

\makeindex

\numberwithin{equation}{section} \numberwithin{theorem}{section}

\begin{document}

\maketitle{\footnotesize We study the conditions under which an
almost Hermitian structure $(G,J)$ of general natural lift type on
the cotangent bundle $T^*M$ of a Riemannian manifold $(M,g)$ is
K\" ahlerian. First, we obtain the algebraic conditions under
which the manifold $(T^*M,G,J)$ is almost Hermitian.  Next we get
the integrability conditions for the almost complex structure $J$,
then the conditions under which the associated $2$-form is closed.
The manifold $(T^*M,G,J)$ is K\" ahlerian iff it is almost
K\"ahlerian and the almost complex structure $J$ is integrable. It
follows that the family of K\"ahlerian structures of above type on
$T^*M$ depends on three essential parameters (one is a certain
proportionality factor, the other two are parameters involved in
the definition of $J$).

\vskip2mm


\par}

\vskip3mm
\section{\bf INTRODUCTION}
\vskip2mm The fundamental differences between the geometry of the
cotangent bundle $T^*M$ and that of the tangent bundle $TM$ of a
Riemannian manifold $(M,g)$, are due to the  different construction
of lifts to $T^*M$, which cannot be defined just like in the case of
$TM$(see [14]).

The possibility to consider vertical, complete and horizontal
lifts on $T^*M$ leads to interesting geometric structures, studied
in the last years (see \cite{Munteanu}, \cite{OprPap1},
\cite{OprPap2}, \cite{OprPap3}, \cite{OprPapMitric},
\cite{OprPor2}, \cite{OprPor1} , \cite{Porosniuc5},
\cite{Porosniuc6}, \cite{Porosniuc3}, \cite{Porosniuc1},
\cite{Porosniuc2}). Due to the duality tangent bundle-cotangent
bundle, some of these results are very much similar to previous
results for the tangent bundles.

In the present paper we study the conditions under which a structure
$(G,J)$ of general natural lift type on the cotangent bundle $T^*M$
of a Riemannian manifold $M$ is a K\"ahlerian structure. The similar
problem for the tangent bundle was treated by Oproiu, in the paper
\cite{Oproiu4} (see also \cite{Oproiu3} for the diagonal case). A
part of the results from the papers \cite{OprPor2},
\cite{Porosniuc6}, \cite{Porosniuc3}, may be obtained as particular
cases from those in the present paper.

First we consider a tensor field $J$ of type $(1,1)$ on $T^*M$
which is a general natural lift of the Riemannian metric $g$. The
condition for $J$ to define an almost complex structure on $T^*M$
leads to certain algebraic relations between the parameters
involved in its  definition. Four of the eight parameters involved
in the definition of $J$ may be obtained as (rational) functions
of the other four parameters. The integrability condition for the
almost complex structure $J$ implies that the base manifold must
have constant sectional curvature. Then some other relations
fulfilled by the parameters and their derivatives are obtained, so
that the essential parameters involved in the definition of the
integrable almost complex structure are two.

In the next section we study the conditions under which a
Riemannian metric $G$ which is a general natural lift of $g$, is
Hermitian with respect to $J$. We get that the structure $(G,J)$
on $T^*M$ is almost Hermitian if and only if the coefficients
$c_1,\  c_2,\ c_3,$ involved in the definition of the metric $G$
are proportional to the coefficients $a_1,\ a_2,\ a_3,$ involved
in definition of  the almost complex structure $J$, with the
proportionality factor $\lambda,$ and the combinations
$c_1+2td_1,\ c_2+2td_2,\ c_3+2td_3$ are proportional to
$a_1+2tb_1,\ a_2+2tb_2,\ a_3+2tb_3$, the proportionality factor
being chosen of the form $\lambda+2t\mu$. The main result is that
the structure $(G,J)$ is K\" ahlerian if and only if
$\mu=\lambda^\prime$ and the coefficients $b_1,\ b_2,\ b_3$ of the
almost complex structure $J$ may be expressed as certain rational
functions of $a_1,\ a_2,\ a_3,$ and their derivatives. The
condition for $(G,J)$ to be K\"ahlerian is obtained from the
property for $(G,J)$ to be almost K\"ahlerian and the
integrability property for $J$.

Some quite long computations have been done by using the Mathematica
package RICCI for doing tensor calculations.

The manifolds, tensor fields and other geometric objects  we
consider in this paper are assumed to be differentiable of class
$C^\infty $ (i.e. smooth). We use the computations in local
coordinates in a fixed local chart but many results may be
expressed in an invariant form by using the vertical and
horizontal lifts. The well known summation convention is used
throughout this paper, the range of the indices $h,i,j,k,l,m,r $
being always $\{1,\dots ,n\}$.

\vskip3mm
\section{\bf PRELIMINARY RESULTS}
\vskip2mm

Let $(M,g)$ be a smooth $n$-dimensional Riemannian manifold and
denote its cotangent bundle by $\pi :T^*M\longrightarrow M$.
Recall that there is a structure of a $2n$-dimensional smooth
manifold on $T^*M$, induced from the structure of smooth
$n$-dimensional manifold  of $M$. From every local chart
$(U,\varphi )=(U,x^1,\dots ,x^n)$  on $M$, it is induced a local
chart $(\pi^{-1}(U),\Phi )=(\pi^{-1}(U),q^1,\dots , q^n,$
$p_1,\dots ,p_n)$, on $T^*M$, as follows. For a cotangent vector
$p\in \pi^{-1}(U)\subset T^*M$, the first $n$ local coordinates
$q^1,\dots ,q^n$ are  the local coordinates of its base point
$x=\pi (p)$ in the local chart $(U,\varphi )$ (in fact we have
$q^i=\pi ^* x^i=x^i\circ \pi, \ i=1,\dots n)$. The last $n$ local
coordinates $p_1,\dots ,p_n$ of $p\in \pi^{-1}(U)$ are the vector
space coordinates of $p$ with respect to the natural basis
$(dx^1_{\pi(p)},\dots , dx^n_{\pi(p)})$, defined by the local
chart $(U,\varphi )$,\ i.e. $p=p_idx^i_{\pi(p)}$.

An $M$-tensor field of type $(r,s)$ on $T^*M$ is defined by sets
of $n^{r+s}$ components (functions depending on $q^i$ and $p_i$),
with $r$ upper indices and $s$ lower indices, assigned to induced
local charts $(\pi^{-1}(U),\Phi )$ on $T^*M$, such that the local
coordinate change rule is that of the local coordinate components
of a tensor field of type $(r,s)$ on the base manifold $M$ (see
\cite{Mok} for further details in the case of the tangent bundle).
An usual tensor field of type $(r,s)$ on $M$ may be thought as an
$M$-tensor field of type $(r,s)$ on $T^*M$. If the considered
tensor field on $M$ is covariant only, the corresponding
$M$-tensor field on $T^*M$ may be identified with the induced
(pullback by $\pi $) tensor field on $T^*M$.

Some useful $M$-tensor fields on $T^*M$ may be obtained as
follows. Let $v,w:[0,\infty ) \longrightarrow {\bf R}$ be smooth
functions and let $\|p\|^2=g^{-1}_{\pi(p)}(p,p)$ be the square of
the norm of the cotangent vector $p\in \pi^{-1}(U)$ ($g^{-1}$ is
the tensor field of type (2,0) having the components
$(g^{kl}(x)),$ which are the entries of the inverse of the matrix
$(g_{ij}(x)),$ defined by the components of $g$ in the local chart
$(U,\varphi )$). The components $vg_{ij}(\pi(p))$, $p_i$,
$w(\|p\|^2)p_ip_j $ define respective $M$-tensor fields of types
$(0,2)$, $(0,1)$, $(0,2)$ on $T^*M$. Similarly, the components
$vg^{kl}(\pi(p))$, $g^{0i}=p_hg^{hi}$, $w(\|p\|^2)g^{0k}g^{0l}$
define respective $M$-tensor fields of type $(2,0)$, $(1,0)$,
$(2,0)$ on $T^*M$. Of course, all the above components are
considered in the induced local chart $(\pi^{-1}(U),\Phi)$.

\vskip3mm
The Levi Civita connection $\dot \nabla $ of $g$ defines
a direct sum decomposition
\begin{equation}\label{splitting}
TT^*M=VT^*M\oplus HT^*M.
\end{equation}
of the tangent bundle to $T^*M$ into vertical distributions
$VT^*M= {\rm Ker}\ \pi _*$ and the horizontal distribution
$HT^*M$.

If $(\pi^{-1}(U),\Phi)=(\pi^{-1}(U),q^1,\dots ,q^n,p_1,\dots
,p_n)$ is a local chart on $T^*M$, induced from the local chart
$(U,\varphi )= (U,x^1,\dots ,x^n)$, the local vector fields
$\frac{\partial}{\partial p_1}, \dots , \frac{\partial}{\partial
p_n}$ on $\pi^{-1}(U)$ define a local frame for $VT^*M$ over $\pi
^{-1}(U)$ and the local vector fields $\frac{\delta}{\delta
q^1},\dots ,\frac{\delta}{\delta q^n}$ define a local frame for
$HT^*M$ over $\pi^{-1}(U)$, where
$$
\frac{\delta}{\delta q^i}=\frac{\partial}{\partial
q^i}+\Gamma^0_{ih} \frac{\partial}{\partial p_h},\ \ \ \Gamma
^0_{ih}=p_k\Gamma ^k_{ih}
 $$
and $\Gamma ^k_{ih}(\pi(p))$ are the Christoffel symbols of $g$.

The set of vector fields $\{\frac{\partial}{\partial p_1},\dots
,\frac{\partial}{\partial p_n}, \frac{\delta}{\delta q^1},\dots
,\frac{\delta}{\delta q^n}\}$ defines a local frame on $T^*M$,
adapted to the direct sum decomposition (\ref{splitting}).

We consider
\begin{equation}
t=\frac{1}{2}\|p\|^2=\frac{1}{2}g^{-1}_{\pi(p)}(p,p)=\frac{1}{2}g^{ik}(x)p_ip_k,
\ \ \ p\in \pi^{-1}(U),
\end{equation}
the energy density defined by $g$ in the cotangent vector $p$. We
have $t\in [0,\infty)$ for all $p\in T^*M$.

Roughly speaking, the natural lifts have coefficients as functions
of the density energy only ([1], [2], [3]).

From now on we shall work in a fixed local chart $(U,\varphi)$ on
$M$ and in the induced local chart $(\pi^{-1}(U),\Phi)$ on $T^*M$.

We may easily obtain the following

\begin{lemma}\label{lema1}
If $n>1$ and $u,v$ are smooth functions on $T^*M$ such that either
$$
u g_{ij}+v p_ip_j=0, \quad  u g^{ij}+v g^{0i}g^{0j}=0,\quad {\rm
or} \quad u\delta ^i_j+vg^{0i} p_j=0
$$
on the domain of any induced local chart on $T^*M$, then $u=0,\
v=0$.
\end{lemma}

\vskip3mm
\section{\bf \ THE INTEGRABILITY  OF THE ALMOST COMPLEX STRUCTURES}
\vskip3mm

If $X$ is a vector field $X\in \mathcal{X}(M)$, then $g_X$ is the
1-form on $M$ defined by the relation $g_X(Y)=g(X,Y),\ \forall Y\in
\mathcal{X}(M)$. Using the local chart of $M$, $(U,x^1,\dots,x^n)$,
we may write $X=X^i\frac{\partial}{\partial x^i}$. Then
$g_X=g_{ij}X^jdx^i,\ g_\frac{\partial}{\partial x^i}=g_{ij}dx^j$.

Let us consider a $1$-form $\theta \in \Lambda^{1}(M)$. Then
$\theta^\sharp=g^{-1}_\theta$ is a vector field on $M$ defined by
the musical isomorphism $g(\theta^\sharp,Y)=\theta (Y), \ \forall
Y\in\mathcal{X}(M)$. If $\theta=\theta_idx^i$, then
$\theta^\sharp= g^{ij}\theta_j \frac{\partial}{\partial x^i}$, the
expressions being on the local chart $(U,x^1,\dots,x^n)$ of $M$.

\vskip2mm

For $p \in T^*M$ we consider the vector $p^\sharp$ tangent to $M$
in $\pi (p)$.

The Liouville vector field on $T^*M$ is $p^V$, expressed by $p_p^V
= p_i \frac{\partial}{\partial p_i} $ in every point $p$ of the
induced local chart $(\pi^{-1}(U),\Phi)$ on $T^*M$.

The similar horizontal vector field on $T^*M$ is $(p^\sharp)^H$,
with the expression $(p^\sharp_p)^H = g^{0i}\frac{\delta}{\delta
q^i}$ in every point $p$ of the induced local chart
$(\pi^{-1}(U),\Phi)$ on $T^*M$ (recall that $g^{0i}=p_hg^{hi}$).

\vskip3mm Consider the real valued smooth functions $a_1,\ a_2,\
a_3,\ a_4,\ b_1,\ b_2,\ b_3,\ b_4$ defined on $[0,\infty)\subset
{\bf R}$. We define \emph{a general natural tensor of type $(1,1)$
on} $T^*M$ by its action on the horizontal and vertical lifts on
$T^*M$
\begin{equation*}
\left\{
\begin{array}{l}
JX^H_p=a_1(t)
(g_X)^V_p+b_1(t)p(X)p_p^V+a_4(t)X_p^H+b_4(t)p(X)(p^\sharp)_p^H,
\\ \mbox{ }  \\
J\theta^V_p=a_3(t)\theta^V_p+b_3(t)g^{-1}_{\pi(p)}
(p,\theta)p_p^V-a_2(t)(\theta^\sharp)_p^H- b_2(t)g^{-1}_{\pi(p)}
(p,\theta)(p^\sharp)_p^H,
\end{array}
\right.
\end{equation*}
in every point $p$ of the induced local card $(\pi^{-1}(U),\Phi)$
on $T^*M$, $\forall ~X \in \mathcal{X}(M), \forall~ \theta \in
\Lambda^1 (M)$.

With respect to the local adapted frame $\{\frac{\delta}{\delta
q^i},\frac{\partial}{\partial p_j}\}_{i,j=1,\dots,n}$, the
expressions of J become

\begin{equation}\label{sist3}
\left\{
\begin{array}{l}
J\frac{\delta}{\delta q^i}=a_1(t)g_{ij}\frac{\partial}{\partial
p_j}+ b_1(t)p_iC+a_4(t)\frac{\delta}{\delta
q^i}+b_4(t)p_i\widetilde C,
\\ \mbox{ }  \\
J\frac{\partial}{\partial p_i}=a_3(t)\frac{\partial}{\partial
p_i}+ b_3(t) g^{0i}C-a_2(t)g^{ij}\frac{\delta}{\delta
q^j}-b_2(t)g^{0i}\widetilde C,
\end{array}
\right.
\end{equation}
where we have denoted by $C =p^V$ the Liouville vector-field on
$T^*M$ and by $\widetilde C=(p^\sharp)^H,$ the corresponding
horizontal vector field on $T^*M$.

We may write also
\begin{equation}\label{sist4}
\left\{
\begin{array}{l}
J\frac{\delta}{\delta q^i}=J^{(1)}_{ij}\frac{\partial}{\partial
p_j}+ J4^j_i\frac{\delta}{\delta q^j},
\\ \mbox{ }  \\
J\frac{\partial}{\partial p_i}=J3^i_j\frac{\partial}{\partial
p_j}-J_{(2)}^{ij}\frac{\delta}{\delta q^j},
\end{array}
\right.
\end{equation}
where
$$
J^{(1)}_{ij}=a_1(t)g_{ij}+b_1(t)p_ip_j,\quad
J4^j_i=a_4(t)\delta^j_i+b_4(t)g^{0j}p_i,
$$
\vskip1mm
$$
\quad J3^i_j=a_3(t)\delta^i_j+b_3(t)g^{0i}p_j,\quad
J_{(2)}^{ij}=a_2(t)g^{ij}+b_2(t)g^{0i}g^{0j}.
$$

\begin{theorem}\label{th4}
A natural tensor field $J$ of type $(1,1)$ on $T^*M$ given by
$(\ref{sist3})$ or $(\ref{sist4})$ defines an almost complex
structure on $T^*M$, if and only if $a_4=-a_3, b_4=-b_3$ and the
coefficients $a_1,\ a_2,\ a_3,\ b_1,\ b_2$ and $b_3$ are related
by
\begin{equation}\label{rel4}
a_1a_2=1+a_3^2\ ,\ \ \ (a_1+2tb_1)(a_2+2tb_2)=1+(a_3+2tb_3)^2.
\end{equation}
\end{theorem}

\vskip 3mm \it{Remark.} \rm From the conditions (\ref{rel4}) we
get that the coefficients $a_1,\ a_2,\  a_1+2tb_1,\ a_2+2tb_2$
have the same sign and cannot vanish. We assume that $a_1>0,\
a_2>0,\ a_1+2tb_1>0,\ a_2+2tb_2>0$ for all $t\geq 0$.

\vskip2mm \it{Remark.} \rm The relations (\ref{rel4}) allow us to
express two of the coefficients $a_1,\ a_2,\ a_3,\ b_1$, $b_2,\
b_3,$ as functions of the other four; e.g. we have:
\begin{equation}\label{inlocuire}
a_2=\frac{1+a_3^2}{a_1},\ \ \
b_2=\frac{2a_3b_3-a_2b_1+2tb_3^2}{a_1+2tb_1}.
\end{equation}

The integrability condition for the above almost complex structure
$J$ on a manifold $M$ is characterized by the vanishing of its
Nijenhuis tensor field $N_J$, defined by
$$
N_J(X,Y)=[JX,JY]-J[JX,Y]-J[X,JY]-[X,Y],
$$
for all vector fields $X$ and $Y$ on $M$.

\begin{theorem}\label{th3}
Let $(M,g)$ be an $n(>2)$-dimensional connected Riemannian
manifold. The almost complex structure $J$ defined by
{\rm(\ref{sist3})} on $T^*M$ is integrable if and only if $(M,g)$
has constant sectional curvature $c$ and the  coefficients $b_1,\
b_2,\ b_3$ are given by:
\begin{equation}\label{integrab}
\begin{cases}
 b_1=\frac{2 c^2 t a_2^2+2 c t a_1 a_2^\prime+a_1
a_1^\prime -c+3 c a_3^2}{a_1-2 t a_1^\prime-2 c t a_2-4 c t^2
a_2^\prime},\\
 \\
b_2=\frac{2 t a_3^{\prime 2}- 2 t a_1^\prime a_2^\prime+c a_2^2+2
c t a_2 a_2^\prime+a_1 a_2^\prime}{a_1-2 t a_1^\prime-2 c t a_2-4
c t^2 a_2^\prime},\\
\\
 b_3 =\frac{a_1 a_ 3^\prime+ 2 c a_2 a_3+4 c t a_2^\prime a_3-
 2 c t a_ 2 a_3^\prime}{a_1-2 t a_1^\prime-2 c t a_2-4 c t^2 a_2^\prime}.
\end{cases}
\end{equation}
\end{theorem}

\it{Proof.} \rm From the condition $N_J(\frac{\partial}{\partial
p_i},\frac{\partial}{\partial p_j})=0$ we obtain that the
horizontal component of this Nijenhuis bracket vanishes if and
only if
\begin{equation}\label{a2p}
a'_2=\frac{a_2a'_3+2a_3b_2-a_2b_3}{2(a_3+tb_3)}
\end{equation}
and the vertical component vanishes if and only if
\begin{equation}\label{curvature}
(a_1a'_2-a_1b_2+2a'_3b_3t)(\delta^h_jp_i-\delta^h_ip_j)-
\end{equation}
$$
a_2^2g^{0k}R^h_{kij}-a_2b_2g^{0k}g^{0l}(p_jR^h_{kil}-
p_iR^h_{kjl})=0.
$$

Taking into account that the curvature of the base manifold does
not depend on $p$, we differentiate with respect to $p_k$ in
(\ref{curvature}). Considering the value of this derivative  in
$p=0$, we get
\begin{equation}\label{sectional curvature}
R^h_{kij}=c(\delta^h_ig_{kj}-\delta^h_jg_{ki}),
\end{equation}
where
$$
c=\frac{a_1(0)}{a_2^2(0)}(b_2(0)-a'_2(0)),
$$
which is a function depending on $q^1,...,q^n$ only. According to
the Schur's theorem, $c$ must be a constant when $n>2$ and $M$ is
connected.

Using the condition of constant sectional curvature for the base
manifold, we obtain from $N_J(\frac{\delta}{\delta
q^i},\frac{\delta}{\delta q^j})=0,$ the following expresions for
$a'_1$ and $a'_3$
\begin{equation}\label{vala1p}
\begin{cases}
a'_1=\frac{a_1b_1+c(1-3a_3^2-4ta_3b_3)}{a_1+2tb_1},\\
\\
a'_3=\frac{a_1b_3-2ca_2(a_3+tb_3)}{a_1+2tb_1}.
\end{cases}
\end{equation}

If we replace in (\ref{a2p}) the expression of $a'_3$ and the
relations (\ref{inlocuire}), we may write $a'_2$ in the form
\begin{equation}\label{vala2p}
a'_2=\frac{2a_3b_3-a_2b_1-ca_2^2}{a_1+2tb_1}.
\end{equation}

The values of $a'_1,\ a'_3$ in (\ref{vala1p}) and $a'_2$ in
(\ref{vala2p}) fulfill the vanishing condition for the vertical
component of the Nijenhuis bracket $N_J(\frac{\partial}{\partial
p_i}, \frac{\delta}{\delta q^j})$. The same expressions fulfill
also the relation
\begin{equation}\label{produs}
a_1a'_2+a'_1a_2=2a_3a'_3,
\end{equation}
obtained by differentiating the first of the relation (\ref{rel4})
with respect to t.

We may solve the system given by (\ref{a2p}) and (\ref{vala1p}),
with respect to $b_1,\ b_2,\ b_3$. Taking (\ref{produs}) into
account, we get the relations (\ref{integrab}) from the theorem,
which fulfill identically the expression of $b_2$ from
(\ref{inlocuire}).

\bigskip
\it{Remark.} \rm In the diagonal case, where $a_3=0$ it follows
$b_3=0$ too, and we have:
$$
a_2=\frac{1}{a_1},\ b_1=\frac{a_1a_1^\prime-c}{a_1-2ta_1^\prime},\
b_2=\frac{c-a_1a_1^\prime}{a_1(a_1^2-2ct)}.
$$

Thus we found some results from \cite{OprPor1} and
\cite{Porosniuc3}, papers which have treated only the diagonal
case.

\section{\bf K\" AHLER  STRUCTURES OF GENERAL NATURAL LIFT TYPE
ON THE COTANGENT BUNDLE}
\vskip3mm  In this section, we introduce
\emph{a Riemannian metric $G$ of general natural lift type on the
cotangent bundle} $T^*M$, defined by
\begin{equation}\label{Ginvar}
\left\{
\begin{array}{l}
G_p(X^H, Y^H) = c_1(t)g_{\pi(p)}(X,Y) + d_1(t)p(X)p(Y),
\\ \mbox{ } \\
G_p(\theta^V,\omega^V) = c_2(t)g^{-1}_{\pi(p)}(\theta,\omega) +
d_2(t)g^{-1}_{\pi(p)}(p,\theta)g^{-1}_{\pi(p)}(p,\omega),
\\ \mbox{ } \\
G_p(X^H,\theta^V) = G_p(\theta^V,X^H)
=c_3(t)\theta(X)+d_3(t)p(X)g^{-1}_{\pi(p)}(p,\theta),
\end{array}
\right.
\end{equation}
$\forall~ X,Y \in \mathcal{X}(M),$ $\forall~ \theta, \omega \in
\Lambda^1(M), \forall~p \in T^*M$.

\vskip3mm

The expressions (\ref{Ginvar}) in local coordinates become
\begin{equation}\label{rel11}
\left\{
\begin{array}{l}
G(\frac{\delta}{\delta q^i}, \frac{\delta}{\delta
q^j})=c_1(t)g_{ij}+ d_1(t)p_ip_j=G^{(1)}_{ij},
\\   \mbox{ } \\
G(\frac{\partial}{\partial p_i}, \frac{\partial}{\partial
p_j})=c_2(t)g^{ij}+d_2(t)g^{0i}g^{0j}=G_{(2)}^{ij},
\\   \mbox{ } \\
G(\frac{\partial}{\partial p_i},\frac{\delta}{\delta q^j})=
G(\frac{\delta}{\delta q^i},\frac{\partial}{\partial p_j})=
c_3(t)\delta_i^j+d_3(t)p_ig^{0j}=G3_i^j,
\end{array}
\right.
\end{equation}
where $c_1,\ c_2,\ c_3,\ d_1,\ d_2,\ d_3$ are six smooth functions
of the density energy on $T^*M$. The conditions for $G$ to be
positive definite are assured if
\begin{equation}\label{pozdef}
c_1+2td_1>0,\quad c_2+2td_2>0,
\end{equation}
$$
(c_1+2td_1)(c_2+2td_2)-(c_3+2td_3)^2>0.
$$

The metric $G$ is almost Hermitian with respect to the general
almost complex structure $J$, if
$$
G(JX,JY)=G(X,Y),
$$
for all vector fields $X,Y$ on $T^*M$.

We  prove the following result

\begin{theorem}\label{th4}{The family of natural Riemannian metrics  $G$
on $T^*M$ such that $(T^*M,G,J)$ is an almost Hermitian manifold,
is given by {\rm(\ref{rel11})}, provided that the coefficients
$c_1,\ c_2,\ c_3,\ d_1,\ d_2,$ and $d_3$ are related to the
coefficients $a_1,\ a_2,\ a_3,\ b_1,\ b_2,$ and $b_3$ by the
following proportionality relations
\begin{equation}\label{proportiec}
\frac{c_1}{a_1} =\frac{c_2}{a_2}=\frac{c_3}{a_3} = \lambda
\end{equation}
\begin{equation}\label{proportied}
\frac{c_1+2 t d_1}{a_1+2 t b_1} =\frac{c_2+2 t d_2}{a_2+2 t
b_2}=\frac{c_3+2 t d_3}{a_3+2 t b_3} = \lambda +2 t \mu,
\end{equation}
where the proportionality coefficients $\lambda>0 $ and $\lambda
+2t\mu>0$ are functions depending on $t$.}
\end{theorem}

\it{Proof.} \rm We use the local adapted frame
$\{\frac{\delta}{\delta q^i},\frac{\partial}{\partial
p_j}\}_{i,j=1,\dots,n}$. The metric $G$ is almost hermitian, if
and only if the following conditions are fulfilled
\begin{equation}\label{aprherm}
\left\{%
\begin{array}{ll}
   G(J\frac{\delta}{\delta q^i},J\frac{\delta}{\delta
q^j})=G(\frac{\delta}{\delta q^i},\frac{\delta}{\delta q^j}), \\
\\
    G(J\frac{\partial}{\partial p_i},J\frac{\partial}{\partial
p_j})=G(\frac{\partial}{\partial p_i},\frac{\partial}{\partial
p_j}),  \\
\\
   G(J\frac{\partial}{\partial p_i},J\frac{\delta}{\delta
q^j})=G(\frac{\partial}{\partial p_i},\frac{\delta}{\delta q^j}).
\end{array}%
\right.
\end{equation}
By using Lemma \ref{lema1} we obtain that the coefficients of
$g_{ij},\ g^{ij},\ \delta^i_j$ in the conditions (\ref{aprherm})
must vanish. It follows that the parameters $c_1,\ c_2,\ c_3,$
from the definition of the metric $G$, satisfy the homogeneous
linear system of the form
\begin{equation}\label{sistc}
\left\{%
\begin{array}{ll}
 (a_3^2-1)c_1+a_1^2c_2-2a_1a_3c_3=0 , \\
 \\
  a_2^2c_1+(a_3^2-1)c_2-2a_2a_3c_3=0 ,  \\
  \\
  a_2a_3c_1+a_1a_3c_2-2a_1a_2c_3=0.
\end{array}%
\right.
\end{equation}

The nontrivial solutions of (\ref{sistc}) are given by the
expression (\ref{proportiec}).

From the vanishing condition for the coefficients of $p_ip_j,\
g^{0i}g^{0j},\ g^{0i}p_j$ in (\ref{aprherm}), we obtain a much
 more complicated system, fulfilled by $d_1,\ d_2,\ d_3$. In order to get
a certain similitude with the above system (\ref{sistc}),
fulfilled by $c_1,\ c_2,\ c_3$, we multiply the new equations by
$2t$ and substract the equations of the system (\ref{sistc}),
respectively. The new system may be written in a form in which the
new unknowns are $c_1+2td_1,\ c_2+2td_2,\ c_3+2td_3:$
\begin{equation*}\label{sistd}
\left\{%
\begin{array}{ll}
 [(a_3+2tb_3)^2-1](c_1+2td_1)+(a_1+2tb_1)^2(c_2+2td_2)\\
 -2(a_1+2tb_1)(a_3+2tb_3)(c_3+2td_3)=0 , \\
 \\
  (a_2+2tb_2)^2(c_1+2td_1)+[(a_3+2tb_3)^2-1](c_2+2td_2)\\
  -2(a_2+2tb_2)(a_3+2tb_3)(c_3+2td_3)=0 ,  \\
  \\
  (a_2+2tb_2)(a_3+2tb_3)(c_1+2td_1)+(a_1+2tb_1)(a_3+2tb_3)(c_2+2td_2)\\
  -2(a_1+2tb_1)(a_2+2tb_2)(c_3+2td_3)=0.  \\
\end{array}%
\right.
\end{equation*}

Then the  nonzero solutions are given by the relation
(\ref{proportied}).

The conditions (\ref{pozdef}) are fullfield, due to the properties
(\ref{inlocuire}) of the coefficients of the almost complex
structure $J$.

The explicit expressions of the coefficients $d_1,d_2,d_3$,
obtained from (\ref{proportied}), are
\begin{equation}\label{proportie2}
\left\{%
\begin{array}{ll}
d_1=\lambda b_1+\mu(a_1+2tb_1), \\
d_2=\lambda b_2+\mu(a_2+2tb_2)  ,  \\
d_3=\lambda b_3+\mu(a_3+2tb_3).
\end{array}%
\right.
\end{equation}

\it{Remark.} \rm In the case when $a_3=0$, it follows that
$c_3=d_3=0$ and we obtain the almost Hermitian structure
considered in \cite{OprPor1} and \cite{Porosniuc3}. Moreover, if
$\lambda =1$ and $\mu =0$, we obtain the almost K\"ahlerian
structure considered in the mentioned papers.

\vskip 3mm Considering the two-form $\Omega $ defined by the
almost Hermitian structure $(G,J)$ on $T^*M$
$$
\Omega (X,Y)=G(X,JY),
$$
for all vector fields $X,Y$ on $T^*M$, we obtain the following
result:

\begin{proposition}\label{prop5}
The expression of the $2$-form $\Omega $ in the local adapted
frame $\{\frac{\delta}{\delta q^i},\frac{\partial}{\partial
p_j}\}_{i,j=1,\dots,n}$ on $T^*M$, is given by
$$
\Omega \left(\frac{\partial}{\partial
p_i},\frac{\partial}{\partial p_j}\right)=0,\ \Omega \left(
\frac{\delta}{\delta q^i},\frac{\delta}{\delta q^j}\right)=0,\
\Omega\left(\frac{\partial}{\partial p_i},\frac{\delta}{\delta
q^j}\right)= \lambda \delta^i_j+\mu g^{0i}p_j
$$
or, equivalently
\begin{equation}
\Omega =(\lambda \delta^i_j+\mu g^{0i}p_j)Dp_i\wedge dq^j,
\end{equation}
where  $Dp_i=dp_i-\Gamma^0_{ih}dq^h$ is the absolute differential
of $p_i$.
\end{proposition}

Next, by calculating the exterior differential of $\Omega$,  we
obtain:
\begin{theorem}\label{th6}
The almost Hermitian structure $(T^*M,G,J)$ is almost
K\"{a}hlerian if and only if
$$
\mu=\lambda ^\prime .
$$
\end{theorem}
\it{Proof.} \rm The differential of $\Omega$ is
$$
d\Omega=(d\lambda \delta^i_j+d\mu g^{0i}p_j+\mu d g^{0i}p_j+\mu
g^{0i}dp_j)\wedge Dp_i\wedge dq^j + (\lambda \delta^i_j+\mu
g^{0i}p_j)dDp_i\wedge dq^j.
$$

We first obtain the expressions of $d\lambda,\ d\mu, \ dg^{0i}$ and
$dDp_i$:
$$
d\lambda=\lambda ' g^{0h}Dp_h,\quad d\mu=\mu ' g^{0h}Dp_h,\quad
dg^{0i}=g^{hi}Dp_h-\Gamma^i_{j0}dq^j,
$$
$$
dDp_i=\frac{1}{2}R^0_{ijh}dq^h\wedge dq^j-\Gamma^h_{ij}Dp_h\wedge
dq^j.
$$

By substituting these relations in $d\Omega$, using the properties
of the external product, the symmetry of $g^{ij},\ \Gamma^h_{ij}$
and the Bianchi identities, we get
$$
d\Omega=\frac{1}{2}(\lambda
'-\mu)p_k(g^{kh}\delta^i_j-g^{ki}\delta^h_j)Dp_h\wedge Dp_i\wedge
dq^j.
$$

Hence $d\Omega=0$ if and only if $\mu=\lambda '$.

\vskip 2mm \it{Remark.} \rm Thus the family of general natural
almost K\"ahlerian structures on $T^*M$ depends on five essential
coefficients $a_1,\ a_3, \ b_1,\ b_3,\ \lambda$, which must
satisfy the supplementary conditions $a_1>0,\ a_1+2tb_1>0,\
\lambda>0, \ \lambda +2 t \mu>0$.

\vskip 3mm Combining the results from the theorems \ref{th4},
\ref{th3} and \ref{th6} we may state

\begin{theorem}
An almost Hermitian structure $(G,J)$ of general natural lift type
on $T^*M$ is K\"ahlerian if and only if  the almost complex
structure $J$ is integrable (see Theorem $\ref{th3}$) and
$\mu=\lambda^\prime$.
\end{theorem}

\it{Remark.} \rm The family of general natural K\"ahlerian
structures on $T^*M$ depends on three  essential coefficients
$a_1,\ a_3,\ \lambda$, which must satisfy the supplementary
conditions $a_1>0,\ a_1+2tb_1>0,\ \lambda>0$,
$\lambda+2t\lambda^\prime>0$, where $b_1$ is given by
(\ref{integrab}).

Examples of such structures may be found in \cite{OprPor2},
\cite{OprPor1}, \cite{Porosniuc3}.

\vskip2mm

\textbf{Acknowledgements.} The author expresses her gratitude to
her PhD adviser, professor V.Oproiu, for several discussions on
this argument, several hints and encouragements.

The research was partially supported by the Grant, TD 158/2007,
CNCSIS, Ministerul Educa\c tiei \c si Cercet\u arii, Rom\^ania.

\begin{flushright}
\it "Al.I. Cuza" University, \\
Faculty of Mathematics, \\
Bd. Carol I, Nr. 11,\\
Ia\c si, 700 506, ROM\^{A}NIA \\
simonadruta@yahoo.com
\end{flushright}

\end{document}